# SCALING LIMIT FOR TRAP MODELS ON $\mathbb{Z}^D$

By Gérard Ben Arous and Jiří Černý

*Courant Institute for Mathematical Sciences and École Polytechnique Fédérale de Lausanne*

We give the "quenched" scaling limit of Bouchaud's trap model in $d \geq 2$. This scaling limit is the fractional-kinetics process, that is the time change of a $d$-dimensional Brownian motion by the inverse of an independent $\alpha$-stable subordinator.

**1. Introduction.** This work establishes scaling limits for certain important models of trapped random walks on $\mathbb{Z}^d$. More precisely we show that Bouchaud's trap model on $\mathbb{Z}^d$, $d \geq 2$, properly normalized, converges (at the process level) to the fractional-kinetics process, which is a self-similar non-Markovian continuous process, obtained as the time change of a $d$-dimensional Brownian motion by the inverse of an independent Lévy $\alpha$-stable subordinator. This is in sharp contrast to the scaling limit for the same model in dimension one (see [8]) where the limiting process is a singular diffusion in random environment. For a general survey about trap models and their motivation in statistical physics we refer to the lecture notes [2], where we announced the result proved in this paper.

Bouchaud's trap model on $\mathbb{Z}^d$ is the nearest neighbor continuous time Markov process $X(t)$ given by the jump rates

$$(1) \qquad c(x, y) = \frac{1}{2d\tau_x} \qquad \text{if } x \text{ and } y \text{ are neighbors in } \mathbb{Z}^d,$$

and zero otherwise, where $\{\tau_x : x \in \mathbb{Z}^d\}$ are i.i.d. heavy-tailed random variables. More precisely we assume that for some $\alpha \in (0, 1)$

$$(2) \qquad \mathbb{P}[\tau_x \geq u] = u^{-\alpha}(1 + L(u)) \qquad \text{with } L(u) \to 0 \text{ as } u \to \infty.$$

We will always assume that $X(0) = 0$. The Markov process $X(t)$ waits at a site $x$ an exponentially distributed time with mean $\tau_x$, and then it jumps to one of the neighbors of $x$ with uniform probability. Therefore $X$ is a random time change of a standard discrete time simple random walk on $\mathbb{Z}^d$. More precisely:

DEFINITION 1.1. Let $S(0) = 0$ and let $S(k)$, $k \in \mathbb{N}$, be the time of the $k$th jump of $X$. For $s \in \mathbb{R}$ we define $S(s) = S(\lfloor s \rfloor)$. We call $S(s)$ the *clock process*. Define the embedded discrete time Markov chain $Y(k)$ by $Y(k) = X(t)$ for $S(k) \leq t < S(k+1)$. Then obviously, $Y$ is a simple random walk on $\mathbb{Z}^d$.

In order to state our principal result we need to introduce the limiting fractional-kinetics (FK) process.

DEFINITION 1.2. Let $B_d(t)$ be a standard $d$-dimensional Brownian motion started at 0, and let $V_\alpha$ be an independent $\alpha$-stable subordinator satisfying $\mathbb{E}[e^{-\lambda V_\alpha(t)}] = e^{-t\lambda^\alpha}$. Define the generalized right-continuous inverse of $V_\alpha(t)$ by $V_\alpha^{-1}(s) := \inf\{t : V_\alpha(t) > s\}$. We define the *fractional-kinetics process* $Z_{d,\alpha}$ by

$$(3) \qquad\qquad Z_{d,\alpha}(s) = B_d(V_\alpha^{-1}(s)).$$

This process is well known in the physics literature. See, for instance, the broad survey by Zaslavsky [20] or the recent book [21] about the relevance of this process for chaotic deterministic systems; see also [9, 10, 13, 14, 17] for more on this class of processes and references.

We fix a time $T > 0$ and $d \geq 2$ and denote by $D^d([0,T])$ the space of càdlàg functions from $[0,T]$ to $\mathbb{R}^d$. Let $X_N(t)$ be the sequence of elements of $D^d([0,T])$,

$$(4) \qquad\qquad X_N(t) = \frac{\sqrt{d}X(tN)}{f(N)},$$

where

$$(5) \qquad f(N) = \begin{cases} C_2(\alpha)N^{\alpha/2}(\log N)^{(1-\alpha)/2}, & \text{if } d = 2, \\ C_d(\alpha)N^{\alpha/2}, & \text{if } d \geq 3, \end{cases}$$

$$(6) \qquad C_d(\alpha) = \begin{cases} [\pi^{1-\alpha}\alpha^{\alpha-1}\Gamma(1-\alpha)\Gamma(1+\alpha)]^{-1/2}, & \text{if } d = 2, \\ [G_d(0)^\alpha\Gamma(1-\alpha)\Gamma(1+\alpha)]^{-1/2}, & \text{if } d \geq 3, \end{cases}$$

and $G_d(0)$ denotes the Green's function of the $d$-dimensional discrete simple random walk at the origin, $G_d(0) = \sum_{k=0}^\infty \mathbb{P}[Y(k) = 0]$, for $d \geq 3$.

Our main result is the following "quenched" scaling limit statement:

THEOREM 1.3. *For a.e.* $\boldsymbol{\tau}$, *the distribution of* $X_N$ *converges weakly to the distribution of* $Z_{d,\alpha}$ *on* $D^d([0,T])$ *equipped with the uniform topology.*

This result is a consequence of the following, more detailed statement, that is, the joint convergence of the clock process and of the position of the embedded random walk. We use $D([0,T], M_1)$ (resp. $D([0,T], U)$) to denote the space $D([0,T])$ equipped with the $M_1$ (resp. uniform) topology. Define

$$(7) \qquad Y_N(t) = \frac{\sqrt{d}}{f(N)} Y(\lfloor tf(N)^2 \rfloor) \quad \text{and} \quad S_N(t) = \frac{1}{N} S(\lfloor tf(N)^2 \rfloor).$$

THEOREM 1.4. *For a.e. $\boldsymbol{\tau}$, the joint distribution of $(S_N, Y_N)$ converges weakly to the distribution of $(V_\alpha, B_d)$ on $D([0,T], M_1) \times D^d([0,T], U)$.*

Let us insist on the following important facts:

1. One word of caution is necessary about the nature of this joint convergence. It takes place in the uniform topology for the spatial component but only in the Skorokhod $M_1$ topology for the clock process (see [18] for the classical reference about the various topologies on $D^d([0,T])$ and [19] for a thorough, more recent, survey). It is important to remark that our statement is not true in the stronger $J_1$ topology (usually called the Skorokhod topology). Indeed, the main advantage of the $M_1$ topology over the $J_1$ topology, for our purposes, is that existence of "intermediate jumps" forbid convergence in the latter but not in the former. These intermediate jumps are important in our context: they are caused by the fact that the deep traps giving the main contributions to the clock process are visited at several nearby instants. All these visits are summed up into one jump of the limiting $\alpha$-stable subordinator $V_\alpha$.

2. Our scaling limit result is "quenched," that is, it is valid almost surely in the random environment $\boldsymbol{\tau}$, and the limiting process is independent of $\boldsymbol{\tau}$.

3. Our result might be seen as a "triviality" result. Indeed, the fractional kinetics process can be obtained as a scaling limit of a much simpler discrete process, that is, a continuous time random walk (CTRW) à la Montroll–Weiss [15]. More precisely consider a simple random walk $Y$ on $\mathbb{Z}^d$ and a sequence of positive i.i.d. random variables $\{s_i : i \in \mathbb{N}\}$ satisfying the same condition (2) as the $\tau_x$'s. Define the CTRW $U(t)$ by

$$(8) \qquad U(t) = Y(k) \qquad \text{if } t \in \left[ \sum_{i=1}^{k-1} s_i, \sum_{i=1}^{k} s_i \right).$$

It is proved in [16] on the level of fixed-time distributions and in [12] on the level of processes that there is a constant $C$ such that

$$(9) \qquad CN^{-\alpha/2} U(tN) \overset{N \to \infty}{\longrightarrow} Z_{d,\alpha}(t).$$

The result of Theorem 1.3 shows that the limit of the $d$-dimensional trap model and its clock process on $\mathbb{Z}^d$ is trivial, in the sense that it is identical with the scaling limit of the much simpler (completely annealed) dynamics of the CTRW. The necessary scaling is the same as for CTRW if $d \geq 3$, and it requires a logarithmic correction if $d = 2$.

4. As mentioned above, the situation is completely different in $d = 1$, where the scaling limit is a singular diffusion in random environment introduced in by Fontes, Isopi and Newman [8] as follows. Let $(x_i, v_i)$ be an inhomogeneous Poisson point process on $\mathbb{R} \times (0, \infty)$ with intensity measure $dx\, \alpha v^{-1-\alpha}\, dv$, and consider the random discrete measure $\rho = \sum_i v_i \delta_{x_i}$ which can be obtained as a scaling limit of the random environment $\boldsymbol{\tau}$. Conditionally on $\rho$, the FIN diffusion $Z_\alpha(s)$ is defined as a diffusion process [with $Z_\alpha(0) = 0$] that can be expressed as a time change of the standard one-dimensional Brownian motion $B_1$ with the speed measure $\rho$: denoting by $\ell(t, y)$ the local time of the standard Brownian motion $B_1$, let

$$(10) \qquad \phi_\rho(t) = \int_{\mathbb{R}} \ell(t, y) \rho(dy),$$

then $Z_\alpha(s) = B(\phi_\rho^{-1}(s))$.

Observe that both processes, the fractional kinetics and the FIN diffusion, are defined as a time change of the Brownian motion $B_d(t)$. The clock processes however differ considerably. For $d = 1$, the clock equals $\phi_\rho(t) = \int \ell(t, y) \rho(dy)$. Therefore, the processes $B_1$ and $\phi_\rho$ are *dependent*. In the fractional-kinetics case the Brownian motion $B_d$ and the clock process, that is, the stable subordinator $V_\alpha$, are *independent*. The asymptotic independence of the clock process $S$ and the location $Y$ is a very remarkable feature distinguishing $d \geq 2$ and $d = 1$.

Note also that, in contrast with the $d = 1$ case, nothing like a scaling limit of the random environment appears in the definition of $Z_{d,\alpha}$ for $d \geq 2$, and that the convergence holds $\boldsymbol{\tau}$-a.s. The absence of the scaling limit of the environment in the definition of $Z_{d,\alpha}$ translates into the non-Markovianity of $Z_{d,\alpha}$. It is, however, considerably easier to control the behavior of the FK process than of the FIN diffusion even if the former is not Markovian. Let us mention few elementary properties of the process $Z_{d,\alpha}$.

PROPOSITION 1.5.  (i) $Z_{d,\alpha}$ *is a.s. $\gamma$-Hölder continuous for any $\gamma < \alpha/2$.*

(ii) $Z_{d,\alpha}$ *is self-similar,* $Z_{d,\alpha}(t) \overset{\text{law}}{=} \lambda^{-\alpha/2} Z_{d,\alpha}(\lambda t)$.

(iii) $Z_{d,\alpha}$ *is not Markovian.*

(iv) *The fixed-time distribution of $Z_{d,\alpha}(t)$ is given by its Fourier transform*

$$(11) \qquad \mathbb{E}(e^{i\xi \cdot Z_{d,\alpha}(t)}) = E_\alpha(-|\xi|^2 t^\alpha / 2),$$

*where $E_\alpha(z) = \sum_{m=0}^{\infty} z^m / \Gamma(1 + m\alpha)$ is the Mittag–Leffler function.*

Proof. Since the Brownian motion is $\gamma$-Hölder continuous for $\gamma < 1/2$ and $V_\alpha^{-1}$ is $\gamma$-Hölder continuous for $\gamma < \alpha$ (see Lemma III.17 of [4]), fact (i) follows. (ii) can be proved using scaling properties of $B_d$ and $V_\alpha$. To show (iii) it is enough to observe that the times between jumps of $Z_{d,\alpha}$ have no exponential distribution. Example B on page 453 of [6] implies that the Laplace transform of $V_\alpha^{-1}(t)$ is equal to $E_\alpha(-\lambda t^\alpha)$. The result of (iv) then follows by an easy computation. $\quad\square$

The name of the FK process comes from the fact that $Z_{d,\alpha}$ has a smooth density $p(t,x)$ which satisfies the fractional-kinetics equation (see [20])

$$(12) \qquad \frac{\partial^\alpha}{\partial t^\alpha} p(t,x) = \frac{1}{2}\Delta p(t,x) + \delta(0)\frac{t^{-\alpha}}{\Gamma(1-\alpha)}.$$

The FK process $Z_{d,\alpha}$ has an obvious aging property due to its very slow clock, namely

$$\mathbb{P}[Z_{d,\alpha}(t_w + s) = Z_{d,\alpha}(t_w) \; \forall s \le t] = \frac{\sin\alpha\pi}{\pi} \int_0^{t_w/(t_w+t)} u^{\alpha-1}(1-u)^{-\alpha}\,du.$$

(13)

This is simply a restatement of the arcsine law for the stable subordinator $V_\alpha$ since

$$(14) \qquad \mathbb{P}[Z_{d,\alpha}(t_w + s) = Z_{d,\alpha}(t_w) \; \forall s \le t] = \mathbb{P}[\{V(t):t\in\mathbb{R}\} \cap [t_w, t_w+t] = \varnothing].$$

This and Theorem 1.3 explain in part the analogous aging result for Bouchaud's trap model

$$(15) \qquad \lim_{t_w\to\infty}\mathbb{P}[X(t_w+\theta t_w) = X(t_w)|\boldsymbol{\tau}] = \frac{\sin\alpha\pi}{\pi}\int_0^{1/(1+\theta)} u^{\alpha-1}(1-u)^{-\alpha}\,du.$$

In fact proving (15) requires a slightly more detailed understanding of the discrete clock process (see [1, 3, 5]).

At the end of the Introduction, we would like to draw reader's attention to the paper [7], where the scaling limit of the trap model on a large complete graph is identified. The situation there is slightly different since there is no natural scaling limit of the simple random walk on a large complete graph in the absence of trapping.

The rest of the paper is organized as follows. In Section 2 we recall the coarse-graining construction introduced for $d=2$ in [3] and we state (for all $d \ge 2$) some results related to this construction. Using these results we prove Theorems 1.3 and 1.4 in Section 3. In Section 4 we give the proofs of the claims from Section 2 for $d \ge 3$.

**2. Coarse graining.** We define in this section the coarse-graining procedure that was used in [3, 5] to prove aging (15). We also recall some properties of this procedure which we need to prove our scaling limit results.

We use $D_x(r)$ [resp. $B_x(r)$] to denote the ball (box) with radius (side) $r$ centered at $x$. These sets are understood as subsets of $\mathbb{Z}^d$. We will often use the claim that $D_x(r)$ contains $d^{-1}\omega_d r^d$ sites, where $\omega_d$ is the surface of a $d$-dimensional unit sphere, although it is not precisely true. Any error we introduce by this consideration is negligible for $r$ large. If $x$ is the origin, we omit it from the notation.

It follows from the definition of $X$ that the clock process $S$ can be written as

$$(16) \qquad S(k) = \sum_{i=0}^{k-1} e_i \tau_{Y_i},$$

where the $e_i$'s are mean-one i.i.d. exponential random variables. We always suppose that the $e_i$'s are coupled with $X$ and $Y$ in this way.

Let $n \in \mathbb{N}$ large. We will consider the processes $Y$ and $X$ before the first exit from the large ball $\mathbb{D}(n) = D(mr(n))$, where (the scales for $d = 2$ are chosen to agree with [3])

$$(17) \qquad r(n) = \begin{cases} \pi^{-1/2} 2^{n/2} n^{(1-\alpha)/2}, & \text{if } d = 2, \\ 2^{n/2}, & \text{if } d \geq 3, \end{cases}$$

and $m$ is a large constant independent of $n$ which will be chosen later (see Corollary 2.2 and Lemma 2.4 below). Let $\zeta_n$ be the exit time of $Y$ from $\mathbb{D}(n)$,

$$(18) \qquad \zeta_n = \inf\{k \in \mathbb{N} : Y(k) \notin \mathbb{D}(n)\}.$$

In $\mathbb{D}(n)$, a principal contribution to the clock process comes from traps with depth of order $g(n)$ where

$$(19) \qquad g(n) = \begin{cases} n^{-1} 2^{n/\alpha}, & \text{if } d = 2, \\ 2^{n/\alpha}, & \text{if } d \geq 3. \end{cases}$$

We define, as in [3],

$$(20) \qquad T_\varepsilon^M(n) = \{x \in \mathbb{D}(n) : \varepsilon g(n) \leq \tau_x < Mg(n)\}.$$

If $M$ or $\varepsilon$ are omitted, it is understood $M = \infty$, respectively $\varepsilon = 0$. We always suppose that $\varepsilon < 1 < M$. We further introduce two $d$-dependent constants $\kappa$, $\gamma$. For $d = 2$ we choose

$$(21) \qquad \gamma < 1 - \alpha \quad \text{and} \quad \kappa = \frac{5}{1-\alpha},$$

for $d \geq 3$

$$(22) \qquad \gamma = 1 - \frac{1}{3d} \quad \text{and} \quad \kappa = \frac{1}{d}.$$

We then define the coarse-graining scale $\rho(n)$ as

$$(23) \qquad \rho(n) = \begin{cases} \pi^{-1/2} 2^{n/2} n^{\gamma/2}, & \text{if } d = 2, \\ 2^{\gamma n/2}, & \text{if } d \geq 3. \end{cases}$$

We will often abbreviate

$$(24) \qquad h(n) = r(n)/\rho(n).$$

The last scale we need is the "proximity" scale

$$(25) \qquad \nu(n) = \begin{cases} \pi^{-1/2} 2^{n/2} n^{-\kappa/2}, & \text{if } d = 2, \\ 2^{\kappa n/2}, & \text{if } d \geq 3. \end{cases}$$

Observe that $\nu(n) \ll \rho(n)$. We use $\mathcal{E}(n)$, $\mathcal{B}(n)$ to denote the sets

$$(26) \quad \mathcal{E}(n) = \{x \in \mathbb{D}(n) : \text{dist}(x, T_\varepsilon^M(n)) > \nu(n)\},$$

$$(27) \quad \mathcal{B}(n) = \begin{cases} \varnothing, & \text{if } d = 2, \\ \{x \in T_\varepsilon^M(n) : (\exists y \neq x : y \in T_\varepsilon^M(n), \text{dist}(x,y) \leq \nu(n))\}, \\ \qquad \text{if } d \geq 3. \end{cases}$$

For all objects defined above we will often skip the dependence on $n$ in the notation.

We now introduce the coarse-graining procedure. Let $j_i^n$ be a sequence of stopping times for $Y$ given by $j_0^n = 0$ and

$$(28) \qquad j_i^n = \min\{k > j_{i-1}^n : \text{dist}(Y(k), Y(j_{i-1}^n)) > \rho(n)\}, \qquad i \in \mathbb{N}.$$

For every $i \in \mathbb{N}_0$ we define the *score* of the part of the trajectory between $j_i^n$ and $j_{i+1}^n$ as follows. Let

$$(29) \qquad \lambda_{i,1}^n = \min\{k \geq j_i^n : Y(k) \in T_\varepsilon^M\}$$

and $y_i^n = Y(\lambda_{i,1}^n)$. Let further

$$(30) \qquad \begin{aligned} \lambda_{i,2}^n &= \min\{k \geq \lambda_{i,1}^n : \text{dist}(Y(k), y_i^n) > \nu(n)\}, \\ \lambda_{i,3}^n &= \min[\{k \geq \lambda_{i,1}^n : Y(k) \in T_\varepsilon^M \setminus y_i^n\} \cup \{k \geq \lambda_{i,2}^n : Y(k) \in T_\varepsilon^M\}]. \end{aligned}$$

If the part of the trajectory between $j_i^n$ and $j_{i+1}^n$ satisfies

$$(31) \qquad \text{dist}(Y(j_i^n), \partial \mathbb{D}(n)) > \rho(n), \qquad Y(j_i^n), Y(j_{i+1}^n) \in \mathcal{E}(n)$$

and

$$(32) \qquad \begin{aligned} \lambda_{i,1}^n &< \lambda_{i,2}^n < j_{i+1}^n \leq \lambda_{i,3}^n, \\ \text{dist}(y_i^n, \partial D_{Y(j_i^n)}(\rho(n))) &> \nu(n), \qquad y_i^n \notin \mathcal{B}(n), \end{aligned}$$

then we define the score of this part as

$$(33) \qquad s_i^n = \sum_{k=\lambda_{i,1}^n}^{\lambda_{i,2}^n} e_k \tau_k \mathbb{1}\{Y(k) = y_i^n\}.$$

If (31) is satisfied and $\lambda_{i,1}^n \geq j_{i+1}^n$, we set $s_i^n = 0$. In both these cases $s_i^n$ records the time spent by $X$ in $T_\varepsilon^M$ during the $i$th part of the trajectory. In all other cases we set $s_i^n = \infty$. This value marks the part of trajectory where something "bad" happens. We use $J(n)$ to denote the index of the first bad part,

$$(34) \qquad\qquad J(n) = \min\{i : s_i^n = \infty\}.$$

We finally introduce two families of random variables, $s^n(x) \in [0, \infty)$ and $r^n(x) \in \mathbb{Z}^d$, indexed by $x \in \mathbb{D}(n)$. By definition, the law of $s^n(x)$ is the same as the law of $s_i^n$ conditioned on $Y(j_i^n) = x$ (and on $\boldsymbol{\tau}$). Similarly, the law of $r^n(x)$ is the same as the law of $Y(j_{i+1}^n) - Y(j_i^n)$ conditioned on the same event.

We will need these properties of the random variables $s^n(x)$.

Lemma 2.1.  *Let*

$$(35) \qquad\qquad \mathcal{E}_0(n) = \{x \in \mathcal{E}(n) : \operatorname{dist}(x, \partial \mathbb{D}(n)) > \rho(n)\}.$$

*Then, for every $\varepsilon$, $M$ and for $\mathbb{P}$-a.e. random environment $\boldsymbol{\tau}$:*

(i)

$$(36) \qquad\qquad \max_{x \in \mathcal{E}_0(n)} \mathbb{P}[s^n(x) = \infty | \boldsymbol{\tau}] = o(h(n)^{-2}).$$

(ii)

$$(37) \quad \lim_{n \to \infty} \max_{x \in \mathcal{E}_0(n)} |h(n)^2 \{1 - \mathbb{E}[e^{-\lambda s^n(x)/2^{n/\alpha}} | s^n(x) < \infty, \boldsymbol{\tau}]\} - F_d(\lambda)| = 0.$$

*Here*

$$(38) \qquad F_d(\lambda) = \mathcal{K}_d \left\{ p_\varepsilon^M - \int_\varepsilon^M \frac{\alpha}{1 + \mathcal{K}_d' \lambda z} \cdot \frac{1}{z^{\alpha+1}} \, dz \right\},$$

$p_\varepsilon^M = \varepsilon^{-\alpha} - M^{-\alpha}$ *and*

$$(39) \qquad \mathcal{K}_d = \begin{cases} (\log 2)^{-1}, \\ 1, \end{cases} \qquad \mathcal{K}_d' = \begin{cases} \pi^{-1} \log 2, & \text{if } d = 2, \\ G_d(0), & \text{if } d \geq 3, \end{cases}$$

(iii)

$$(40) \qquad\qquad \lim_{n \to \infty} \max_{x \in \mathcal{E}_0(x)} |h(n)^2 \mathbb{P}[s^n(x) \neq 0 | \boldsymbol{\tau}] - \mathcal{K}_d p_\varepsilon^M| = 0.$$

For $d = 2$ (i) follows from Section 5, (ii) from Lemma 6.4 and (iii) from Lemma 5.7 of [3]. We give in Section 4 a proof for $d \geq 3$ taken from [5].

It is worth noting that (i) of the previous lemma implies that (ii) holds also when conditioning on $s^n(x) < \infty$ is removed. As a corollary of (i) we also get

COROLLARY 2.2. *For every $\delta, T > 0$ there exists $m$ independent of $\varepsilon$ and $M$, such that $\boldsymbol{\tau}$-a.s. for $n$ large*

$$\mathbb{P}[J(n)/h(n)^2 \geq T | \boldsymbol{\tau}] \geq 1 - \delta. \tag{41}$$

PROOF. By (2)

$$\mathbb{P}[0 \notin \mathcal{E}_0(n)] \leq \sum_{x \in D(\nu(n))} \mathbb{P}[x \in T_\varepsilon^M] \leq C\nu(n)^d g(n)^{-\alpha}. \tag{42}$$

This is $O(n^{-\kappa+\alpha})$ for $d = 2$ and $O(2^{-n/2})$ for $d \geq 3$. In both cases the Borel–Cantelli lemma implies that $\boldsymbol{\tau}$-a.s. $0 \in \mathcal{E}_0(n)$ for $n$ large. Therefore, by Lemma 2.1, $\mathbb{P}[s_0^n = \infty | \boldsymbol{\tau}] = o(h^{-2})$. Moreover, by the second condition in (31), if $s_0^n < \infty$, then the first part of the trajectory ends in $\mathcal{E}$. Actually, it ends in $\mathcal{E}_0$ since the set $\mathcal{E} \setminus \mathcal{E}_0$ is at distance $r - \rho \gg \rho$ from the origin. Therefore, the second part of the trajectory starts in $\mathcal{E}_0$ and thus $s_1^n = \infty$ with probability $o(h^{-2})$. This remains true for all parts of the trajectory before the walk approaches the boundary of $\mathbb{D}$. However, since $r/\rho = h$, the expected number of parts needed to approach $\partial\mathbb{D}$ scales as $m^2 h^2$. Therefore, it is possible to choose $m$ large enough such that $Th^2$ parts stay in $D(mr - \rho)$ with probability larger than $1 - \delta/2$. The probability that at least one of these parts is bad is $Th^2 o(h^{-2}) = o(1)$. This completes the proof. $\square$

The behavior of the random variables $r^n(x)$ is easy to control.

LEMMA 2.3. *For every $\xi \in \mathbb{R}^d$ and for all $x \in \mathcal{E}(n)$*

$$\lim_{n \to \infty} h(n)^2 \{1 - \mathbb{E}(e^{-\xi \cdot r^n(x)/r(n)})\} = -\frac{|\xi|^2}{2d}. \tag{43}$$

PROOF. By definition $|r^n(x)| = \rho(1 + o(1)) = r/h \ll r$. Using the Taylor expansion and the symmetry of the distribution of $r^n(x)$ we get

$$\mathbb{E}[e^{-\xi \cdot r^n(x)/r(n)}] = 1 + \mathbb{E}\left[\frac{1}{2}h(n)^{-2}\left(\xi \cdot \frac{r^n(x)}{\rho(n)}\right)^2\right] + O(h(n)^{-4}). \tag{44}$$

It follows, for example, from Lemma 1.7.4 of [11] that the distribution of $r^n(x)/\rho$ converges to the uniform distribution on the sphere of radius one. The result then follows by an easy integration. $\square$

The reason why the scores $s_i^n$ were introduced in [3] is that the sum of scores is a good approximation for the clock process.

LEMMA 2.4. *For any $\delta > 0$ and $T > 0$ one can choose $\varepsilon$, $M$ and $m$ such that $\boldsymbol{\tau}$-a.s. for all $n$ large enough,*

$$\mathbb{P}\left[\frac{1}{2^{n/\alpha}}\max\left\{\left|S(j_k^n) - \sum_{j=0}^{k-1} s_j^n\right| : k \in \{1, \ldots, h^2 T\}\right\} \geq \delta \Big| \boldsymbol{\tau}\right] < \delta. \tag{45}$$

The proof of this lemma for $d \geq 2$ can be found on pages 30–31 of [3]. For $d \geq 3$ it is proved in Section 4.

**3. Proofs of Theorems 1.3 and 1.4.** We prove Theorem 1.4 first. The next lemma gives the convergence of fixed-time marginals.

LEMMA 3.1. *The finite-dimensional distributions of the pair* $(S_N, Y_N)$ *converge to those of* $(V_\alpha, B_d)$.

In order to prove Lemma 3.1 we will need an important lemma describing the asymptotic behavior of the *joint* Laplace transform of $r^n(x)$ and $s^n(x)$.

LEMMA 3.2. *For* $\mathbb{P}$*-a.e. random environment* $\boldsymbol{\tau}$ *and for all* $\lambda > 0, \xi \in \mathbb{R}^d$

$$
\begin{aligned}
\lim_{n \to \infty} h(n)^2 \Big\{ 1 - \mathbb{E}\Big[ \exp\Big( &-\frac{\lambda s^n(x)}{2^{n/\alpha}} - \frac{\xi \cdot r^n(x)}{r(n)} \Big) \Big| s^n(x) < \infty, \boldsymbol{\tau} \Big] \Big\} \\
&= F_d(\lambda) - \frac{|\xi|^2}{2d}
\end{aligned}
$$

*uniformly in* $x \in \mathcal{E}_0(n)$.

PROOF. Note first that by Lemma 2.1(i) $\mathbb{P}[s^n(x) = \infty | \boldsymbol{\tau}] = o(h(n)^{-2})$. Therefore, we can remove the conditioning on $s^n(x) < \infty$. To shorten the expressions we do not explicitly write conditioning on $\boldsymbol{\tau}$ in this proof. By a trivial decomposition according to the value of $s^n(x)$ we get

$$
\begin{aligned}
\mathbb{E}\Big[ \exp\Big( &-\frac{\lambda s^n(x)}{2^{n/\alpha}} - \frac{\xi \cdot r^n(x)}{r(n)} \Big) \Big] \\
&= \mathbb{E}\Big[ \exp\Big( -\frac{\xi \cdot r^n(x)}{r(n)} \Big) \mathbb{1}\{ s^n(x) = 0 \} \Big] \\
&\quad + \mathbb{E}\Big[ \exp\Big( -\frac{\lambda s^n(x)}{2^{n/\alpha}} \Big) \mathbb{1}\{ s^n(x) \neq 0 \} \Big] \cdot \mathcal{R}(n),
\end{aligned}
$$

where, since $|r^n(x)| = \rho(1 + o(1))$,

$$
e^{-\rho(n)|\xi|/r(n)} \leq \mathcal{R}(n) \leq e^{\rho(n)|\xi|/r(n)}
$$

and therefore $\mathcal{R}(n) = 1 + o(1)$. The first expectation on the right-hand side of (47) can be rewritten using Lemma 2.3,

$$
\begin{aligned}
\mathbb{E}\Big[ \exp\Big( &-\frac{\xi \cdot r^n(x)}{r(n)} \Big) \mathbb{1}\{ s^n(x) = 0 \} \Big] \\
&= \mathbb{E}\Big[ \exp\Big( -\frac{\xi \cdot r^n(x)}{r(n)} \Big) \Big] - \mathbb{E}\Big[ \exp\Big( -\frac{\xi \cdot r^n(x)}{r(n)} \Big) \mathbb{1}\{ s^n(x) \neq 0 \} \Big] \\
&= 1 + \frac{|\xi|}{2dh(n)^2} + o(h(n)^{-2}) - \mathcal{R}(n)\mathbb{P}[s^n(x) \neq 0],
\end{aligned}
$$

where $\mathcal{R}(n)$ satisfies again (48). We rewrite the second expectation of (47) using Lemma 2.1(ii),

(50)
$$\mathbb{E}\left[\exp\left(-\frac{\lambda s^n(x)}{2^{n/\alpha}}\right)\mathbb{1}\{s^n(x)\neq 0\}\right]$$
$$=\mathbb{E}\left[\exp\left(-\frac{\lambda s^n(x)}{2^{n/\alpha}}\right)\right]-\mathbb{E}\left[\exp\left(-\frac{\lambda s^n(x)}{2^{n/\alpha}}\right)\mathbb{1}\{s^n(x)=0\}\right]$$
$$=1-h(n)^{-2}F_d(\lambda)+o(h(n)^{-2})-\mathbb{P}[s^n(x)=0]$$
$$=-h(n)^{-2}F_d(\lambda)+o(h(n)^{-2})+\mathbb{P}[s^n(x)\neq 0].$$

Putting everything together we get

(51)
$$\mathbb{E}\left[\exp\left(-\frac{\lambda s^n(x)}{2^{n/\alpha}}-\frac{\xi\cdot r^n(x)}{\rho(n)h(n)}\right)\right]$$
$$=1+\frac{|\xi|}{2dh(n)^2}-\frac{F_d(\lambda)}{h(n)^2}+o(h(n)^{-2})+\{1-\mathcal{R}(n)\}\mathbb{P}[s^n(x)\neq 0].$$

Since $1-\mathcal{R}(n)=o(1)$ and by Lemma 2.1(iii) $\mathbb{P}[s^n(x)\neq 0]=O(h(n)^{-2})$, the proof is complete. $\square$

PROOF OF LEMMA 3.1. To check the convergence of the finite-dimensional distributions of $(S_N,Y_N)$ we choose $n=n(N)\in\mathbb{N}$ and $t=t(N)\in[1,2^{1/\alpha})$ such that

(52)
$$N=2^{n(N)/\alpha}t(N).$$

It is easy to see from the definitions of $n$, $t$ and $r(n)$ that

(53)
$$f(N)=c_1 r(n(N))t(N)^{\alpha/2},$$

where

(54)
$$c_1=c_1(d,\alpha)=\begin{cases}\pi^{1/2}(\alpha^{-1}\log 2)^{(1-\alpha)/2}, & \text{if }d=2,\\ 1, & \text{if }d\geq 3.\end{cases}$$

We further set $c_2=c_2(d,\alpha)=(C_d(\alpha)c_1(d,\alpha))^{-1}$.

Later we will take the limit $n\to\infty$ for a fixed value of $t\in[1,2^{1/\alpha})$ instead of taking the limit $N\to\infty$. We will show that this limit exists and does not depend on $t$. Moreover, as can be seen from the proof, the convergence is uniform in $t$. Therefore also the limit as $N\to\infty$ exists. We will not comment on the issue of uniformity during the proof. Hence, instead of the convergence of $(S_N,Y_N)$ we show that (in the sense of the finite-dimensional distributions) for all $t\in[1,2^{1/\alpha})$

(55)
$$\left(\frac{1}{t2^{n/\alpha}}S(c_2^{-2}r(n)^2t^\alpha\cdot),\frac{c_2\sqrt{d}}{r(n)t^{\alpha/2}}Y(c_2^{-2}r(n)^2t^\alpha\cdot)\right)\xrightarrow{n\to\infty}(V_\alpha(\cdot),B_d(\cdot)).$$

Let $r_k^n = Y(j_{k+1}^n) - Y(j_k^n)$. We will approximate the processes on the left-hand side of the last display by sum of scores $s_j^n$ and of displacement $r_j^n$. It follows from the properties of the simple random walk that the exit time $j_1^n$ from the ball $D(\rho(n))$ satisfies $\mathbb{E}[j_1^n] = \rho(n)^2(1 + o(1))$ and $\mathbb{E}[(j_1^n/\rho(n)^2)^2] < C$ for some $C$ independent of $n$. Therefore, by the law of large numbers for triangular arrays, a.s. for any $\delta' > 0$, $u \leq T$ and $n$ large enough

$$(56) \qquad j_{\lfloor(1-\delta')c_2^{-2}h(n)^2 t^\alpha u\rfloor}^n \leq c_2^{-2} r(n)^2 t^\alpha u \leq j_{\lfloor(1+\delta')c_2^{-2}h(n)^2 t^\alpha u\rfloor}^n.$$

Since $S(\cdot)$ is increasing, $S(c_2^{-2}r(n)^2 t^\alpha u)$ can be approximated from above and below by $S(j_{\lfloor(1\pm\delta')c_2^{-2}h(n)^2 t^\alpha u\rfloor}^n)$. Lemma 2.4 then yields that for $\varepsilon$ small and $M$, $m$, $n$ large

$$(57) \qquad \mathbb{P}\left[\left|\frac{1}{t^{2n/a}}S(j_{\lfloor(1\pm\delta')c_2^{-2}h(n)^2 t^\alpha u\rfloor}^n) - \sum_{i=0}^{\lfloor(1\pm\delta')c_2^{-2}h(n)^2 t^\alpha u\rfloor-1} s_i^n\right| \geq \delta \,\Big|\, \boldsymbol{\tau}\right] \leq \delta.$$

Similarly, it follows from the properties of the simple random walk that for any $\delta > 0$ it is possible to choose $\delta'$ such that

$$\mathbb{P}\left[\left|\frac{c_2\sqrt{d}}{r(n)t^{\alpha/2}}Y(c_2^{-2}r(n)^2 t^\alpha u) - \frac{c_2\sqrt{d}}{r(n)t^{\alpha/2}}\sum_{i=0}^{\lfloor(1\pm\delta')c_2^{-2}h(n)^2 t^\alpha u\rfloor-1} r_i^n\right| \geq \delta\right] \leq \delta.$$
$$(58)$$

Let $0 = u_0 < u_1 < \cdots < u_q \leq T$, $\lambda_i > 0$ and $\xi_i \in \mathbb{R}^d$, where $i \in \{1, \ldots, q\}$. To prove the convergence of the finite-dimensional distributions we will prove the $\boldsymbol{\tau}$-a.s. convergence of the Laplace transform

$$\mathbb{E}\Bigg[\exp\bigg(-\sum_{i=1}^q \frac{\lambda_i}{t^{2n/\alpha}}\{S(c_2^{-2}r(n)^2 t^\alpha u_i) - S(c_2^{-2}r(n)^2 t^\alpha u_{i-1})\}$$

$$(59)$$
$$+ \frac{c_2\sqrt{d}}{r(n)t^{\alpha/2}}\xi_i \cdot \{Y(c_2^{-2}r(n)^2 t^\alpha u_i) - Y(c_2^{-2}r(n)^2 t^\alpha u_{i-1})\}\bigg)\,\Big|\, \boldsymbol{\tau}\Bigg].$$

The discussion of the last paragraph implies that it suffices to show the convergence of

$$(60) \qquad \mathbb{E}\Bigg[\exp\bigg(-\sum_{i=1}^q \sum_{k\in B_v(n,i)} \frac{\lambda_i}{t^{2n/\alpha}}s_k^n + \frac{c_2\sqrt{d}}{r(n)t^{\alpha/2}}\xi_i \cdot r_k^n\bigg)\,\Big|\, \boldsymbol{\tau}\Bigg]$$

for $v = \pm\delta'$, where

$$(61) \qquad B_v(n,i) = \{\lfloor(1+v)c_2^{-2}h(n)^2 t^\alpha u_{i-1}\rfloor, \ldots, \lfloor(1+v)c_2^{-2}h(n)^2 t^\alpha u_i\rfloor - 1\},$$

and to show that as $\delta' \to 0$ both limits coincide.

Let $\mathcal{Q}_n$ be the set of all finite sequences

$$(62) \qquad \mathcal{Q}_n = \{x_\ell \in \mathbb{Z}^d : \ell \in 0, \ldots, \lfloor c_2^{-2}h(n)^2 t^\alpha u_q\rfloor - 1\}.$$

Expression (60) can be written as

$$\mathbb{E}\left[\exp\left(-\sum_{i=1}^{q}\sum_{k\in B_v(n,i)}\frac{\lambda_i}{t2^{n/\alpha}}s_k^n + \frac{c_2\sqrt{d}}{r(n)t^{\alpha/2}}\xi_i \cdot r_k^n\right)\Big|\boldsymbol{\tau}\right]$$

(63)
$$= \sum_{\{x_\ell\}\in Q_n}\mathbb{P}[Y(j_\ell^n)=x_\ell \ \forall\ell]$$

$$\times\mathbb{E}\left[\exp\left(-\sum_{i=1}^{q}\sum_{k\in B_v(n,i)}\frac{\lambda_i s_k^n}{t2^{n/\alpha}} + \frac{c_2\sqrt{d}}{r(n)t^{\alpha/2}}\xi_i \cdot r_k^n\right)\Big|\boldsymbol{\tau},\right.$$

$$\left.Y(j_\ell^n)=x_\ell \ \forall\ell\right].$$

The last summation can be further divided into two parts. We first consider sequences such that $\{x_\ell\}\not\subset\mathcal{E}_0$. It follows from Lemma 2.1 and Corollary 2.2 that the sum of probabilities of such sequences can be made arbitrarily small by choosing $\varepsilon$, $M$ and $m$. We can therefore ignore them. The contribution of the remaining sequences $\{x_\ell\}$ can be evaluated using Lemma 3.2. Indeed, let $\omega=\lfloor c_2^{-2}h(n)^2t^\alpha u_q\rfloor-1$. Observe that given $\boldsymbol{\tau}$ and $Y(j_\omega^n)=x_\omega$, the distribution of $(s_\omega^n,r_\omega^n)$ is independent of the history of the walk and is the same as the distribution of $(s^n(x_\omega),r^n(x_\omega))$. Therefore,

$$\mathbb{E}\left[\exp\left(-\sum_{i=1}^{q}\sum_{k\in B_v(n,i)}\frac{\lambda_i s_k^n}{t2^{n/\alpha}} + \frac{c_2\sqrt{d}}{r(n)t^{\alpha/2}}\xi_i \cdot r_k^n\right)\Big|\boldsymbol{\tau},Y(j_\ell^n)=x_\ell \ \forall\ell\leq\omega\right]$$

$$= \mathbb{E}\left[\exp\left(-\sum_{i=1}^{q}\sum_{\substack{k\in B_v(n,i)\\k\leq\omega-1}}\frac{\lambda_i s_k^n}{t2^{n/\alpha}} + \frac{c_2\sqrt{d}}{r(n)t^{\alpha/2}}\xi_i \cdot r_k^n\right)\Big|\boldsymbol{\tau},\right.$$

(64)
$$\left.Y(j_\ell^n)=x_\ell \ \forall\ell\leq\omega-1\right]$$

$$\times\mathbb{E}\left[\exp\left(-\frac{\lambda_q}{t2^{n/\alpha}}s^n(x_\omega) - \frac{c_2\sqrt{d}}{r(n)t^{\alpha/2}}\xi_q \cdot r^n(x_\omega)\right)\Big|\boldsymbol{\tau}\right].$$

The last expectation is bounded uniformly in $x_\omega$ by

(65)
$$1-(1\pm\delta)h(n)^{-2}\left(F_d\left(\frac{\lambda_q}{t}\right) - \frac{|\xi_q|^2}{2d}\frac{dc_2^2}{t^\alpha}\right).$$

Therefore, we can sum over $x_\omega$ and repeat the same manipulation for $x_{\omega-1}$. Iterating, we find that the sum over $\{x_\ell\}\subset\mathcal{E}_0$ is bounded from above by

$$\mathbb{P}[Y(j_\ell^n)\in\mathcal{E}_0 \ \forall\ell\leq\omega]$$

$$
\begin{aligned}
(66) \quad & \times \prod_{i=1}^{q}\left\{1+\frac{1+\delta}{h(n)^2}\frac{|\xi_i|^2}{2d}\frac{dc_2^2}{t^\alpha}-\frac{1-\delta}{h(n)^2}F_d\left(\frac{\lambda_i}{t}\right)\right\}^{|B_v(n,i)|} \\
& = \prod_{i=1}^{q}\exp\left\{\left((1+\delta)\frac{|\xi_i|^2}{2}-(1-\delta)\frac{t^\alpha}{c_2^2}F_d\left(\frac{\lambda_i}{t}\right)\right)(u_i-u_{i-1})(1+v)\right\} \\
& \quad \times (1+o(1)).
\end{aligned}
$$

A lower bound can be constructed analogously. Obviously, as $\delta,\delta'\to 0$ and $v=\pm\delta'$ the upper and lower bound coincide.

Finally, taking $\varepsilon=0$ and $M=\infty$ in the definition (38) of $F_d(\lambda)$ we find by an easy integration that

$$
(67) \qquad F_d(\lambda)\xrightarrow[\varepsilon\to 0]{M\to\infty}\mathcal{K}_d(\mathcal{K}_d'\lambda)^\alpha\Gamma(1+\alpha)\Gamma(1-\alpha).
$$

Therefore, using definitions (54) and (6),

$$
(68) \qquad \frac{t^\alpha}{c_2^2}F_d\left(\frac{\lambda}{t}\right)\xrightarrow[\varepsilon\to 0]{M\to\infty}\lambda^\alpha.
$$

This completes the proof of Lemma 3.1.   □

To complete the proof of Theorem 1.4 we need to show the tightness.

LEMMA 3.3. *The sequence of the distributions of $(S_N,Y_N)$ is tight in $D([0,T],M_1)\times D^d([0,T],U)$.*

PROOF. To check the tightness for $S_N$ in $D([0,T],M_1)$ we use Theorem 12.12.3 of [19]. Since the $S_N$ are increasing, it is easy to see that condition (i) of this theorem is equivalent to the tightness of $S_N(T)$ which can be easily checked from Lemma 3.1. In order to check condition (ii) of the theorem remark that for increasing functions the oscillation function $w_s$ used in [19] is equal to zero. So checking (ii) boils down to controlling the boundary oscilations $\bar{v}(x,0,\delta)$ and $\bar{v}(x,T,\delta)$. For the first quantity (using again the monotonicity of $S_N$) this amounts to check that for any $\varepsilon,\eta>0$ there is $\delta$ such that $\mathbb{P}[S_N(\delta)\geq\eta]<\varepsilon$ which follows again from Lemma 3.1. The reasoning for $\bar{v}(x,T,\delta)$ is analogous.

For $Y_n$ the proof of the tightness is analogous to the same proof for Donsker's invariance principle. The tightness of both components implies the tightness of the pair $(S_N,Y_N)$ in the product topology on $D([0,T],M_1)\times D^d([0,T],U)$.   □

Obviously, Lemmas 3.1 and 3.3 imply Theorem 1.4. We can now easily derive Theorem 1.3.

PROOF OF THEOREM 1.3. It is easy to check from definitions (4) and (7) that $X_N(\cdot) = Y_N(S_N^{-1}(\cdot))$. Let $D_{u,\uparrow}$ denote the subset of $D([0, T])$ consisting of unbounded increasing functions. By Corollary 13.6.4 of [19] the inverse map from $D_{u,\uparrow}(M_1)$ to $D_{u,\uparrow}(U)$ is continuous at strictly increasing functions. Since the Lévy process $V_\alpha$ [the limit of $S_N$ in $(D_{u,\uparrow}, M_1)$] is a.s. strictly increasing, the distribution of $S_N^{-1}$ converges to the distribution of $V_\alpha^{-1}$ weakly on $D_{u,\uparrow}(U)$ and the limit is a.s. continuous. The composition $(f, g) \mapsto f \circ g$ as the mapping from $D^d([0, T], U) \times D_{u,\uparrow}(U)$ to $D^d([0, T], U)$ is continuous at $C^d \times C$ (here $C$ is the space of continuous function) as is easy to check. The weak convergence of $X_N$ on $D^d([0, T], U)$ then follows. $\square$

## 4. Proofs of the coarse-graining estimates for $d \geq 3$.
We give here the proofs of Lemmas 2.1 and 2.4 for $d \geq 3$. These proofs are adapted from [5] and use similar techniques as in [3] for $d = 2$. In general, the proofs become slightly simpler because the random walk is transient if $d \geq 3$, and all important quantities (like Green's function, hitting probabilities, etc.) depend only polynomially on the radius of the ball, logarithmic corrections are not required.

4.1. *Proof of Lemma 2.1 for $d \geq 3$.* Lemma 2.1 controls the distribution of the random scores $s^n(x)$ for $x \in \mathcal{E}_0$. Typically, $s^n(x)$ is equal to the time that $X$ started at $x$ spends in $T_\varepsilon^M$ before exiting $D_x(\rho)$. In some exceptional cases $s^n(x) = \infty$. We first show that the probability that this happens is $o(h(n)^{-2})$, that is we prove (i) of Lemma 2.1.

As follows from the definition of $s^n(x)$ the exceptional cases that we need to control are:

(a) the exit point from $D_x(\rho(n))$ is not in $\mathcal{E}(n)$,

(b) $Y$ hits a trap in $T_\varepsilon^M(n)$ that is at distance smaller than $\nu(n)$ from $\partial D_x(\rho(n))$,

(c) $Y$ hits two different traps in $T_\varepsilon^M(n)$ before the exit of $D_x(\rho(n))$,

(d) $Y$ hits a trap from $\mathcal{B}(n)$ [see (27) for definition],

(e) $Y$ hits a trap $y$ in $T_\varepsilon^M(n)$, exits $D_y(\nu(n))$ and then returns to $y$ before exiting $D_x(\rho(n))$.

We now bound the probability of all these events. For the event (a) we have

LEMMA 4.1. *Let $P_1(n, x)$ be the probability that the simple random walk started at $x$ exits $D_x(\rho(n))$ at some site that is not in $\mathcal{E}$. Then $\boldsymbol{\tau}$-a.s. for every $x \in \mathcal{E}_0$, $P_1(n, x) \leq Cg(n)^{-\alpha}\nu(n)^d = o(h(n)^{-2})$.*

PROOF. Let $A_x = A_x(n)$ denote the annulus

$$(69) \qquad A_x(n) = D_x(\rho(n) + \nu(n)) \setminus D_x(\rho(n) - \nu(n)).$$

We first show that there exists $K$ such that $\boldsymbol{\tau}$-a.s. for $n$ large enough

$$(70) \qquad |A_x \cap T_\varepsilon^M| \leq K\rho^{d-1}\nu g^{-\alpha} \qquad \text{for all } x \in \mathbb{D}.$$

The number of the sites in $A_x$ is bounded by $|A_x| \leq c'\rho^{d-1}\nu$. Hence, for $x$ fixed

$$\begin{aligned}
& \mathbb{P}[|A_x \cap T_\varepsilon^M| \geq K\rho^{d-1}\nu g^{-\alpha}] \\
(71) \qquad & \leq \exp(-\lambda K\rho^{d-1}\nu g^{-\alpha})\{1 + cg^{-\alpha}\varepsilon^{-\alpha}(e^\lambda - 1)\}^{c'\rho^{d-1}\nu} \\
& \leq \exp\{\rho^{d-1}g^{-\alpha}\nu[-\lambda K + c(e^\lambda - 1)]\}.
\end{aligned}$$

Summing over $x \in \mathbb{D}$ we bound the probability that (70) is violated by

$$(72) \qquad cr(n)^2 \exp\{\rho^{d-1}g^{-\alpha}\nu[-\lambda K + c(e^\lambda - 1)]\}.$$

Since $\rho^{d-1}g^{-\alpha}\nu = 2^{(d-1)\gamma/2 + n\kappa/2 - 1}$ and $(d-1)\gamma/2 + \kappa/2 - 1 > 0$ for our choice of constants, the fact (70) follows by choosing $K$ large and using the Borel–Cantelli lemma.

If (70) is true, then there are at most $cK\rho^{d-1}\nu g^{-\alpha}\nu^{d-1}$ points on the boundary of $D_x(\rho)$ that are not in $\mathcal{E}$. The probability that $Y$ exits $D_x(\rho)$ in any such point is $O(\rho^{1-d})$ (see [11], Lemma 1.7.4). Hence,

$$(73) \qquad P_1(n, x) \leq cK\rho^{d-1}g^{-\alpha}\nu^d\rho^{1-d} = Cg^{-\alpha}\nu^d = o(h(n)^{-2}).$$

This completes the proof. $\quad\square$

Next, we bound the probability that (b) happens.

LEMMA 4.2. *Let $P_2(n, x)$ be the probability that the simple random walk started at $x$ hits a trap in $T_\varepsilon^M(n) \cap A_x(n)$ before exiting $D_x(\rho(n))$. Then $\boldsymbol{\tau}$-a.s. for all $n$ large, $P_2(n, x) \leq C\rho(n)\nu(n)g(n)^{-\alpha} = o(h(n)^{-2})$ for all $x \in \mathbb{D}$.*

PROOF. According to (70) there are $\boldsymbol{\tau}$-a.s. at most $K\rho^{d-1}\nu g^{-\alpha}$ traps in $A_x \cap T_\varepsilon^M$. The probability that the walk hits one such trap $y$ is by (148) bounded from above by $c|x - y|^{2-d}$. There exists constant $C$ such that for all $y \in A_x$, $|x - y|^{2-d} \leq C\rho^{2-d}$. The required probability is thus smaller than $C\rho^{d-1}\nu g^{-\alpha}\rho^{2-d} \leq C\rho\nu g^{-\alpha}$. $\quad\square$

Let $\mathbb{P}_x$ denote the distribution of the simple random walk $Y$ started from $x$. To proof (c) we need several technical lemmas first.

LEMMA 4.3. *Let*

$$(74) \qquad V_x(n) = \sum_{y \in T_\varepsilon^M} \mathbb{P}_x[Y \text{ hits } y \text{ before exiting } D_x(\rho(n)) | \boldsymbol{\tau}].$$

*Then for any $\delta > 0$ and $\boldsymbol{\tau}$-a.s. there is $n_0$ such that for all $n \geq n_0$ and for all $x \in \mathcal{E}_0(n)$*

$$(75) \qquad (1 - \delta)\mathcal{K}_d p_\varepsilon^M h(n)^{-2} \leq V_x(n) \leq (1 + \delta)\mathcal{K}_d p_\varepsilon^M h(n)^{-2}.$$

Proof. Let $\mu = 1 - 2/(3d)$ and $\iota(n) = 2^{\mu n/2}$, therefore $\kappa < \mu < \gamma$ and $\nu(n) \ll \iota(n) \ll \rho(n)$. Recall that $B_x(r)$ denotes the cube centered at $x$ with side $r$. Let $\mathcal{D}_n = D_x(\rho - 2\iota) \setminus B_x(\iota)$. We divide the sum in (74) into three parts. We use $\Sigma_1$ to denote the sum over $y \in T_\varepsilon^M \cap \mathcal{D}_n$, $\Sigma_2$ to denote the sum over $y \in T_\varepsilon^M \cap B_x(\iota)$, and $\Sigma_3$ to denote the sum over $y \in T_\varepsilon^M \cap (D_x(\rho) \setminus D_x(\rho - 2\iota))$. The reason why we introduce the third sum is the error term in (149) which is too large for the traps that are too close to the border of $D(\rho)$.

The main contribution comes from $\Sigma_1$, so we treat it first. We cover $\mathcal{D}_n$ by cubes with side $\iota$. It is not difficult to show that $\tau$-a.s. for $n$ large

$$(76) \qquad |B_x(\iota) \cap T_\varepsilon^M| \in ((1-\delta)p_\varepsilon^M \iota^d g^{-\alpha}, (1+\delta)p_\varepsilon^M \iota^d g^{-\alpha})$$

for all $x \in \mathbb{D}$ such that $\operatorname{dist}(x, \partial \mathbb{D}) \geq \iota\sqrt{2}$. Indeed, let

$$(77) \qquad F_x = \{|B_x(\iota) \cap T_\varepsilon^M| \geq (1+\delta)p_\varepsilon^M \iota^d g^{-\alpha}\}.$$

Then for any small $\eta$ and $n$ large enough $\mathbb{P}[x \in T_\varepsilon^M] \leq (1+\eta)p_\varepsilon^M g^{-\alpha}$. Hence, for $\lambda > 0$

$$(78) \qquad \begin{aligned} \mathbb{P}[F_x] &\leq \exp(-\lambda(1+\delta)p_\varepsilon^M \iota^d g^{-\alpha})\{1 + (e^\lambda - 1)(1+\eta)p_\varepsilon^M g^{-\alpha}\}^{\iota^d} \\ &\leq \exp\{p_\varepsilon^M \iota^d g^{-\alpha}[-\lambda(1+\delta) + (e^\lambda - 1)(1+\eta)]\}. \end{aligned}$$

For any $\delta$ one can choose $\lambda$ and $\eta$ small enough such that the exponent in the last expression is negative. Hence,

$$(79) \qquad \mathbb{P}[F_x] \leq \exp(-c\iota^d g^{-\alpha})$$

for $n$ large enough. Summing over all $x$ and using the definitions of $\iota$ and $g$ we get

$$(80) \qquad \mathbb{P}\left[\bigcup_x F_x\right] \leq Cr^d \exp(-c2^{n(d\mu - 2)/2}).$$

Since $d\mu - 2 > 0$, the upper bound for (76) is finished. The proof of the lower bound is completely analogous.

We can now actually estimate $\Sigma_1$. Without loss of generality we set $x = 0$. Let

$$(81) \qquad H_n = \{z \in \iota(n)\mathbb{Z}^d \setminus \{0\} : B_z(\iota) \cap \mathcal{D}_n \neq \varnothing\}.$$

Using the bound (149) we get

$$(82) \qquad \begin{aligned} \Sigma_1 &\leq \sum_{y \in T_\varepsilon^M \cap \mathcal{D}_n} a_d\{|y|^{2-d} - \rho^{2-d} + O(|y|^{1-d})\}(1 + O(\rho - |y|)^{2-d}) \\ &\leq \sum_{z \in H_n} \sum_{\substack{y \in T_\varepsilon^M \\ y \in B_z(\iota)}} a_d\{|y|^{2-d} - \rho^{2-d} + O(|y|^{1-d})\}(1 + O(\rho - |y|)^{2-d}), \end{aligned}$$

where $a_d = \frac{d}{2}\Gamma(\frac{d}{2}-1)\pi^{-d/2}$. Obviously, for any $y \in B_z(\iota)$, $||y|^{2-d} - |z|^{2-d}| \le c\iota|z|^{1-d}$. This together with (76) yields the bound

$$
(83) \qquad \Sigma_1 \le \sum_{z \in H_n} (1+\delta)p_\varepsilon^M \iota^d g^{-\alpha} a_d \{|z|^{2-d} - \rho^{2-d} + O(\iota|z|^{1-d})\} + \mathcal{R},
$$

where

$$
(84) \qquad \mathcal{R} = \sum_{z \in H_n} \sum_{\substack{y \in T_\varepsilon^M \\ y \in B_z(\iota)}} a_d \{|y|^{2-d} - \rho^{2-d} + O(|y|^{1-d})\} O(\rho - |y|)^{2-d}.
$$

Every site $y$ from the last summation satisfies $|y| \le \rho - \iota$. Therefore, $O(\rho - |y|)^{2-d} = O(\iota^{2-d})$. The error term $\mathcal{R}$ is thus much smaller than the sum in (83) which we now estimate. Replacing the summation by integration (making again an error of order $\iota|z|^{1-d}$) we get

$$
\Sigma_1 \le (1+\delta)p_\varepsilon^M g^{-\alpha} \int_{\mathcal{D}} a_d \{|z|^{2-d} - \rho^{2-d} + O(\iota|z|^{1-d})\} \, dz + \mathcal{R}
$$

$$
(85) \qquad \le (1+\delta)p_\varepsilon^M g^{-\alpha} \rho^2 a_d \omega_d \left(\frac{1}{2} - \frac{1}{d}\right)(1+o(1))
$$

$$
\le (1+2\delta)\mathcal{K}_d p_\varepsilon^M h(n)^{-2},
$$

where $\omega_d$ denotes as before the surface of the $d$-dimensional unit sphere. The lower bound for $\Sigma_1$ can be obtained in the same way. It is actually much simpler, because the lower bound (147) on the hitting probability is less complicated than the upper bound (149). Hence,

$$
(86) \qquad \Sigma_1 \ge (1-2\delta)\mathcal{K}_d p_\varepsilon^M h(n)^{-2}.
$$

It remains to show that $\Sigma_2$ and $\Sigma_3$ are $o(h(n)^{-2})$. To estimate $\Sigma_2$ we need a finer description of the homogeneity of the environment than (76). Let $i_{\max}$ be the smallest integer satisfying $2^i \nu(n) \ge \iota(n)$, that is, $i_{\max} \sim (\mu - \kappa)n/2$. Then $\boldsymbol{\tau}$-a.s. for $n$ large, all $i \in \{0, \ldots, i_{\max}\}$, and all $x \in \mathbb{D}$

$$
(87) \qquad |B_x(2^i \nu) \cap T_\varepsilon^M| \le n(1 \vee 2^{id} \nu^d g^{-\alpha}).
$$

Indeed, fix $i \in \{-1, \ldots, i_{\max}\}$ first. Then for any $x \in \mathbb{D}$ we have

$$
\mathbb{P}[|B_x(2^{n\gamma+i}) \cap T_\varepsilon^M| \ge n(1 \vee 2^{id} \nu^d g^{-\alpha})]
$$

$$
(88) \qquad \le \exp(-\lambda n(1 \vee 2^{id} \nu^d g^{-\alpha}))\{1 + c(e^\lambda - 1)\varepsilon^{-\alpha} g^{-\alpha}\}^{2^{id} \nu^d}
$$

$$
\le C \exp(-c\lambda n).
$$

Summing over $x \in \mathbb{D}$ and $i \in \{-1, \ldots, i_{\max}\}$ we get an upper bound for the probability of the complement of (87) which is of order $nr(n)^d e^{-\lambda n}$. Therefore, choosing $\lambda$ large enough, (87) is true $\mathbb{P}$-a.s. for $n$ large enough.

Let $E = \{-1, 0, 1\}^d \setminus \{0, 0, 0\}$. Let $\mathcal{O}_i$ be the union of $3^d - 1$ cubes of size $2^i \nu$ centered at $2^i \nu E$,

$$(89) \qquad \mathcal{O}_i = \bigcup_{x \in E} B_{x 2^i \nu}(2^i \nu).$$

To bound $\Sigma_2$ we cover the cube $B(\iota)$ (we suppose again that $x = 0$) by $\bigcup_{i=0}^{i_{\max}} \mathcal{O}_i$. Observe that our covering does not contain $B(\nu)$. However, $B(\nu) \subset D(\nu)$ and $0 \in \mathcal{E}_0$, so that $B(\nu) \cap T_\varepsilon^M = \varnothing$.

By (148) and (87) we get

$$(90) \qquad \begin{aligned} \Sigma_2 &\leq C \sum_{i=0}^{i_{\max}} n(1 \vee 2^{id} \nu^d g^{-\alpha})(2^i \nu)^{(2-d)} \\ &\leq C \sum_{i=0}^{(\mu - \kappa)n/2} n\{(2^i \nu)^{2-d} \vee 2^{2i} \nu^2 g^{-\alpha}\}. \end{aligned}$$

The first term in the braces is decreasing in $i$ and the second one is increasing. The sum is thus bounded by $Cn^2(\nu^{2-d} \vee \iota^2 g^{-\alpha})$. However, both terms, $n^2 \nu^{2-d}$ and $n^2 \iota^2 g^{-\alpha}$, are much smaller than $h(n)^{-2}$ for our choice of constants. This means that $\Sigma_2 \ll \Sigma_1$.

The sum $\Sigma_3$, that is the sum over $y \in T_\varepsilon^M \cap (D(\rho) \setminus D(\rho - 2\iota))$ can be bounded in the same way as the probability of hitting a trap in the annulus $A_x \cap T_\varepsilon^M$ was bounded in Lemma 4.2. Following the same reasoning [with $\nu(n)$ replaced by $\iota(n)$] we get $\Sigma_3 \leq \rho \iota g^{-\alpha} \ll h(n)^{-2}$. This completes the proof of Lemma 4.3. $\square$

The second technical lemma that we need to bound the event (c) also provides the required bound for the event (d).

LEMMA 4.4. *Let*

$$(91) \qquad W_x(n) = \sum_{y \in \mathcal{B}(n)} \mathbb{P}_x[Y \text{ hits } y \text{ before exiting } D_x(\rho(n)) | \boldsymbol{\tau}].$$

*Then $\boldsymbol{\tau}$-a.s., for and for all $x \in \mathcal{E}(n)$ and $n$ large enough $W_x(n) = o(h(n)^{-2})$.*

PROOF. The proof is very similar to the previous one. We divide the sum into three parts in the same way as before. We keep the notation $\Sigma_1$, $\Sigma_3$, $\Sigma_3$ for these parts. Since $\mathcal{B} \subset T_\varepsilon^M$, it follows from the previous proof that $\Sigma_2$ and $\Sigma_3$ are $o(h(n)^{-2})$. Hence, it remains to bound $\Sigma_1$ from above. This can be achieved by the same calculation as before if we show that

$$(92) \qquad |B_x(\iota) \cap \mathcal{B}| = o(|B_x(\iota) \cap T_\varepsilon^M|) = o(\iota^d g^{-\alpha})$$

for all $x \in \mathbb{D}$ [cf. this with (76)]. We will show that $\boldsymbol{\tau}$-a.s. for $n$ large and for all $x \in \mathbb{D}$

$$(93) \qquad |B_x(\iota) \cap \mathcal{B}(n)| \leq n^2 \nu^d g^{-2\alpha} r^d =: \phi(n).$$

This bound is not optimal but sufficient for our purposes. Indeed, using the definitions of $g$, $\nu$ and $r$ we find that $\phi(n) = O(n^2 2^{(d-3)n/2})$ which is much smaller than $\iota^d g^{-\alpha} = O(2^{(d-8/3)n/2})$.

Let $\mathcal{L}_n$ denote the grid $\iota(n)\mathbb{Z}^d$. Then, $|\mathcal{L}_n \cap \mathbb{D}| \leq c(r/\iota)^d$. We use $A$ to denote the event that there exists a cube of side $\iota$ containing more than $\phi(n)$ bad sites. If $A$ is true, then there is also a cube of side $2\iota$ centred on $\mathcal{L}_n$ that contains more than $\phi(n)$ bad sites. Therefore,

$$(94) \qquad \mathbb{P}[A] \leq \sum_{x \in \mathcal{L}_n \cap \mathbb{D}} \mathbb{P}[|B_x(2\iota) \cap \mathcal{B}| \geq \phi(n)] \leq C(r/\iota)^d \mathbb{P}[|B(2\iota) \cap \mathcal{B}| \geq \phi(n)].$$

Using the definition of $\mathcal{B}$ and the union bound we get that $\mathbb{P}[x \in \mathcal{B}] \leq c\varepsilon^{-2\alpha}\nu^d g^{-2\alpha}$. Therefore, by the Markov inequality,

$$(95) \qquad \begin{aligned} \mathbb{P}[|B(2\iota) \cap \mathcal{B}| \geq \phi(n)] &\leq \phi(n)^{-1}\mathbb{E}\left[\sum_{x \in B(2\iota)} \mathbb{1}\{x \in \mathcal{B}\}\right] \\ &\leq C\varepsilon^{-2\alpha} n^{-2}(r/\iota)^{-d}. \end{aligned}$$

Putting this into (94) we obtain $\mathbb{P}[A] \leq Cn^{-2}$. Therefore, (93) follows by the Borel–Cantelli lemma and the proof is complete. $\quad\square$

We now use the last two lemmas to bound the probability of the event (c).

LEMMA 4.5. *Let $P_3(n, x)$ denote the probability that the simple random walk started at $x$ hits two traps from $T_\varepsilon^M(n)$ before exiting $D_x(\rho(n))$. Then $\boldsymbol{\tau}$-a.s. for every $x \in \mathcal{E}_0$, $P_3(n, x) = o(h(n)^{-2})$.*

PROOF. For $A \subset \mathbb{Z}^d$ we use $\mathcal{Y}(x, \rho, A)$ to denote the number of different traps from $A$ visited by the simple random walk $Y$ before the exit from $D_x(\rho)$. Then,

$$(96) \qquad \begin{aligned} P_3(n, x) &= \mathbb{P}_x[\mathcal{Y}(x, \rho, T_\varepsilon^M) \geq 2|\boldsymbol{\tau}] \\ &\leq \mathbb{P}_x[\mathcal{Y}(x, \rho, T_\varepsilon^M) \geq 2|\mathcal{Y}(x, \rho, T_\varepsilon^M \setminus \mathcal{B}) \geq 1, \boldsymbol{\tau}] \\ &\quad \times \mathbb{P}_x[\mathcal{Y}(x, \rho, T_\varepsilon^M \setminus \mathcal{B}) \geq 1|\boldsymbol{\tau}] \\ &\quad + \mathbb{P}_x[\mathcal{Y}(x, \rho, T_\varepsilon^M) \geq 2|\mathcal{Y}(x, \rho, \mathcal{B}) \geq 1, \boldsymbol{\tau}]\mathbb{P}_x[\mathcal{Y}(x, \rho, \mathcal{B}) \geq 1|\boldsymbol{\tau}]. \end{aligned}$$

By Lemma 4.3,

$$(97) \qquad \mathbb{P}_x[\mathcal{Y}(x, \rho, T_\varepsilon^M \setminus \mathcal{B}) \geq 1|\boldsymbol{\tau}] \leq V_x(n) = O(h(n)^{-2})$$

and, by Lemma 4.4,

$$(98) \qquad \mathbb{P}_x[\mathcal{Y}(x,\rho,\mathcal{B}) \geq 1 | \boldsymbol{\tau}] \leq W_x(n) = o(h(n)^{-2}).$$

If we show that

$$(99) \qquad \mathbb{P}_x[\mathcal{Y}(x,\rho,T_\varepsilon^M) \geq 2 | \mathcal{Y}(x,\rho,T_\varepsilon^M \setminus \mathcal{B}) \geq 1, \boldsymbol{\tau}] = O(h(n)^{-2}) = o(1),$$

then the lemma follows from (96)–(99). To prove (99) we denote by $y$ the first visited trap from $T_\varepsilon^M$. Then from the strong Markov property and from $D(x,\rho) \subset D(y,2\rho)$ if follows that

$$(100) \qquad \begin{aligned} &\mathbb{P}_x[\mathcal{Y}(x,\rho,T_\varepsilon^M) \geq 2 | \mathcal{Y}(x,\rho,T_\varepsilon^M \setminus \mathcal{B}) \geq 1, \boldsymbol{\tau}] \\ &\quad \leq \mathbb{P}_y[\mathcal{Y}(y,2\rho,T_\varepsilon^M \setminus \{y\}) \geq 1 | \boldsymbol{\tau}] \\ &\quad \leq \sum_{z \in T_\varepsilon^M} \mathbb{P}_y[Y \text{ hits } z \text{ before exiting } D_y(2\rho(n)) | \boldsymbol{\tau}]. \end{aligned}$$

The right-hand side of the last formula can be bounded by $Ch(n)^{-2}$ using the same argument as in Lemma 4.3. This argument works because $y \in T_\varepsilon^M \setminus \mathcal{B}$ and therefore $(T_\varepsilon^M \cap D_y(\nu)) \setminus \{y\} = \varnothing$. The fact that the ball considered in (100) is two times larger than in Lemma 4.3 does not change the asymptotic behavior, it only changes the prefactor. $\quad\square$

It remains to exclude the event (e).

LEMMA 4.6.  *Let $P_5(n,x)$ denote the probability that the simple random walk started at $x$ hits a trap $y \in T_\varepsilon^M(n)$ exits $D_y(\nu(n))$ and returns to $y$ before exiting $D_x(\rho(n))$. Then $\boldsymbol{\tau}$-a.s. for every $x \in \mathcal{E}_0(n)$, $P_5(n,x) = o(h(n)^{-2})$.*

PROOF.  Due to Lemma 4.2 we can suppose that $\text{dist}(y, \partial D_x(\rho)) \geq \nu$. Let $p_{\text{return}}(x,y)$ denote the probability that the simple random started at $y$ that have exited $D_y(\nu)$ returns to $y$ before exiting $D_x(\rho)$. Obviously $P_5(n,x) \leq \max\{p_{\text{return}}(x,y) : y \in D_x(\rho) \cap T_\varepsilon^M\}$. Let $G_A(x,y)$ denote the Green's function of $Y$ killed on the first exit from the set $A \subset \mathbb{Z}^d$. By the decomposition on the first exit from $D_y(\nu)$,

$$(101) \qquad G_{D_x(\rho)}(y,y) = G_{D_y(\nu)}(y,y) + p_{\text{return}}(x,y) G_{D_x(\rho)}(y,y).$$

Hence, by (145), uniformly for $y \in D_x(\rho) \cap T_\varepsilon^M$ and $\text{dist}(y, \partial D_x(\rho)) \geq \nu$,

$$(102) \qquad \begin{aligned} p_{\text{return}}(x,y) &= 1 - \frac{G_{D_y(\nu)}(y,y)}{G_{D_x(\rho)}(y,y)} \\ &\leq 1 - \frac{G_{D(\nu)}(0,0)}{G_{D(2\rho)}(0,0)} = O(\nu^{2-d}) = o(h(n)^{-2}). \end{aligned}$$

This completes the proof of (e) and therefore also of Lemma 2.1(i) for $d \geq 3$.
□

We now show Lemma 2.1(ii). Since $s^n(x)$ records the time that $X$ spends in $T_\varepsilon^M$, we should first control the distribution of the depth of the first hit trap in $T_\varepsilon^M$. To this end we define

$$(103) \qquad \sigma(n) = n^{-1} + \left( \max_{x \geq \varepsilon} |L(g(n)x)| \right)^{1/2},$$

with $L$ defined in (2). Since $L(x) \to 0$ as $x \to \infty$, the function $\sigma$ satisfies

$$(104) \qquad \sigma(n) \geq 1/n, \qquad \lim_{n \to \infty} \sigma(n) = 0$$

and

$$(105) \qquad \max_{x \geq \varepsilon} |L(g(n)x)| \ll \sigma(n) \qquad \text{as } n \to \infty.$$

Further, let $z_n(i)$ be a sequence satisfying $\varepsilon = z_n(0) < z_n(1) < \cdots < z_n(R_n) = M$, and $z_n(i+1) - z_n(i) \in (\sigma(n), 2\sigma(n))$ for all $i \in \{0, \ldots, R_n - 1\}$. Let $p_i^n$ denote the factor

$$(106) \qquad p_i^n = \frac{1}{z_n(i)^\alpha} - \frac{1}{z_n(i+1)^\alpha}.$$

LEMMA 4.7. *Let $\mathcal{P}_x(n,i)$ denote the probability that the simple random walk started at $x$ hits the set $T_{z_n(i)}^{z_n(i+1)}(n)$ before exiting $D_x(\rho(n))$. Then for any $\delta$ and a.e. $\boldsymbol{\tau}$ there is $n_0$ such that for all $n \geq n_0$, for all $x \in \mathcal{E}_0(n)$ and for all $i \in \{0, \ldots, R_n\}$*

$$(107) \qquad \mathcal{P}_x(n,i) \in (\mathcal{K}_d(1-\delta)h(n)^{-2}p_i^n, \mathcal{K}_d(1+\delta)h(n)^{-2}p_i^n).$$

PROOF. We will need the following technical claims.

LEMMA 4.8. *For any $\delta > 0$, $\boldsymbol{\tau}$-a.s. for all $n$ large:*

(i) $\mathbb{P}[0 \in T_{z_n(i)}^{z_n(i+1)}] \in ((1-\delta)g^{-\alpha}p_i^n, (1+\delta)g^{-\alpha}p_i^n)$,

(ii) *for all $x \in \mathbb{D}$ and $i \in \{0, \ldots, R-1\}$*

$$(108) \qquad |B_x(\iota) \cap T_{z_n(i)}^{z_n(i+1)}| \in ((1-\delta)\iota^d g^{-\alpha}p_i^n, (1+\delta)\iota^d g^{-\alpha}p_i^n).$$

PROOF. By (2) we have

$$(109) \qquad p_i^n = g^{-\alpha} \left[ \left( \frac{1}{z_n(i)^\alpha} - \frac{1}{z_n(i+1)^\alpha} \right) + \frac{L(gz_n(i))}{z_n(i)^\alpha} - \frac{L(gz_n(i+1))}{z_n(i+1)^\alpha} \right].$$

To prove (i) we should thus show that

$$(110) \qquad \frac{L(gz_n(i))}{z_n(i)^\alpha} - \frac{L(gz_n(i+1))}{z_n(i+1)^\alpha} = o\left(\frac{1}{z_n(i)^\alpha} - \frac{1}{z_n(i+1)^\alpha}\right).$$

However, this is true since $z_n(i)^{-\alpha} - z_n(i+1)^{-\alpha} \asymp \sigma(n)$ and, as follows from (105), $L(gz_n(i)) = o(\sigma(n))$.

The claim (ii) can be proved exactly as (76) was proved; the estimate on $\mathbb{P}[x \in T_\varepsilon^M]$ should be replaced by the first claim of the lemma. The easy proof is left to the reader. $\quad\square$

We can now finish the proof of Lemma 4.7. We use $V_{x,i}(n)$ to denote

$$(111) \qquad V_{x,i}(n) = \sum_{y \in T_{z_n(i)}^{z_n(i+1)}} \mathbb{P}_x[Y \text{ hits } y \text{ before exiting } D_x(\rho)|\boldsymbol{\tau}].$$

Lemma 4.8(ii) and the procedure used to show Lemma 4.3 give

$$(112) \qquad (1-\delta)\mathcal{K}_d p_i^n h(n)^{-2} \leq V_{x,i}(n) \leq (1+\delta)\mathcal{K}_d p_i^n h(n)^{-2}.$$

For $\mathcal{P}_x(n,i)$ we have then

$$(113) \qquad \mathcal{P}_x(n,i) \leq V_{x,i}(n) \leq (1+\delta)\mathcal{K}_d p_i^n h(n)^{-2}.$$

The corresponding lower bound can be obtained using Lemma 4.5 and Bonferroni's inequality. Indeed, using the notation introduced before (96),

$$\mathcal{P}_x(n,i) \geq V_{x,i}(n) - \sum_{y,z \in T_{z_n(i)}^{z_n(i+1)}} \mathbb{P}_x[Y \text{ hits } y \text{ and } z \text{ before exiting } D_x(\rho)|\boldsymbol{\tau}]$$

$$(114) \qquad = V_{x,i}(n) - \mathbb{P}[\mathcal{Y}(x,\rho,T_{z_n(i)}^{z_n(i+1)}) \geq 2]$$

$$\geq V_{x,i}(n) - \mathbb{P}[\mathcal{Y}(x,\rho,T_\varepsilon^M) \geq 2] \geq (1-2\delta)\mathcal{K}_d p_i^n h(n)^{-2}.$$

This completes the proof. $\quad\square$

We can now show Lemma 2.1(ii), that is, to show that $\boldsymbol{\tau}$-a.s.

$$(115) \qquad \lim_{n\to\infty} \max_{x \in \mathcal{E}_0(n)} |h(n)^2\{1 - \mathbb{E}[e^{-\lambda s^n(x)/2^{n/\alpha}}|s^n(x) < \infty, \boldsymbol{\tau}]\} - F_d(\lambda)| = 0.$$

When the simple random walk $Y$ hits a deep trap $y$ before exiting $D_x(\rho(n))$ and $s^n(x) < \infty$, then $s^n(x)$ is simply the time spent in $y$ before the exit from $D_y(\nu(n))$. The process $Y$ hits $y$ a geometrical number of times. The mean of this geometrical variable is $G_{D(\nu(n))}(0,0)$. Each visit takes an exponential time with mean $\tau_y$. Using the expression (145) from the Appendix we get the following formula for the conditional Laplace transform of $s^n(x)$:

$$(116) \qquad \mathbb{E}\left[\exp\left(-\frac{\lambda s^n(x)}{2^{n/\alpha}}\right)\Big|\tau_y, s^n(x) < \infty\right] = \frac{1}{1 + \lambda\tau_y 2^{-n/\alpha}G_d(0)(1+o(1))}.$$

The probability that $s^n(x) = \infty$ is $o(h(n)^{-2})$. Therefore,

$$
\begin{aligned}
(117) \quad & \mathbb{E}\left[\exp\left(-\frac{\lambda s^n(x)}{2^{n/\alpha}}\right)\Big| s^n(x) < \infty, \boldsymbol{\tau}\right] \\
& \qquad = \mathbb{E}\left[\exp\left(-\frac{\lambda s^n(x)}{2^{n/\alpha}}\right)\Big| \boldsymbol{\tau}\right](1 + o(h(n)^{-2})).
\end{aligned}
$$

The last expectation can be estimated using Lemma 4.7 and (116),

$$
\begin{aligned}
(118) \quad & \mathbb{E}\left[\exp\left(-\frac{\lambda s^n(x)}{2^{n/\alpha}}\right)\Big| \boldsymbol{\tau}\right] \\
& \qquad \geq (1 - (1+\delta)\mathcal{K}_d p_\varepsilon^M h(n)^{-2}) \\
& \qquad\quad + \mathcal{K}_d h(n)^{-2} \sum_{i=1}^{R_n} \frac{p_i^n(1-\delta)}{1 + \lambda z_n(i)G_d(0)(1+o(1))}.
\end{aligned}
$$

For $n$ large the last expression is bounded from below by

$$
\begin{aligned}
(119) \quad & 1 - \mathcal{K}_d h(n)^{-2}\left(p_\varepsilon^M - \int_\varepsilon^M \frac{\alpha}{1 + \lambda G_d(0)z} \cdot \frac{1}{z^{\alpha+1}}\,dz\right) - \delta C h(n)^{-2} p_\varepsilon^M \\
& \qquad = 1 - h(n)^{-2}(F_d(\lambda) + O(\delta)).
\end{aligned}
$$

This and (117) give an upper bound for $1 - \mathbb{E}[e^{-\lambda s^n(x)/2^{n/\alpha}} | s^n(x) < \infty, \boldsymbol{\tau}]$. A corresponding lower bound can be constructed analogously. This completes the proof of Lemma 2.1(ii).

To prove Lemma 2.1(iii) define first $\mathcal{P}_x(n)$ as the probability that the simple random walk started at $x$ hits the set $T_\varepsilon^M$ before exiting $D_x(\rho(n))$. Using Lemma 4.3 and the same reasoning as in (114) it can be proved that

$$
(120) \quad \mathcal{P}_x(n) \in (\mathcal{K}_d(1-\delta)h(n)^{-2}p_\varepsilon^M, \mathcal{K}_d(1+\delta)h(n)^{-2}p_\varepsilon^M).
$$

Since $\mathbb{P}[s^n(x) \neq 0 | \boldsymbol{\tau}]$ is bounded from below by $\mathcal{P}_x(n)$ and from above by $\mathcal{P}_x(n) + \mathbb{P}[s^n(x) = \infty | \boldsymbol{\tau}]$, Lemma 2.1(iii) follows from Lemma 2.1(i) and (120).

4.2. *Proof of Lemma* 2.4 *for* $d \geq 3$. We want to show that for any $\delta > 0$ and $T > 0$ it is possible to choose $\varepsilon$, $M$ and $m$ such that for $\boldsymbol{\tau}$-a.s. and $n$ large enough

$$
(121) \quad \mathbb{P}\left[\frac{1}{2^{n/\alpha}}\max\left\{\left|S(j_k^n) - \sum_{j=0}^{k-1} s_j^n\right| : k \in \{1, \dots, h(n)^2 T\}\right\} \geq \delta\right] < \delta.
$$

The sum of scores records (if $s_i^n$ stay finite) only the time spent in $T_\varepsilon^M$. Let $\mathcal{G}_n$ be the event $\{s_j^n < \infty : j \leq T h(n)^2\}$. As follows from Corollary 2.2, the probability of $\mathcal{G}_n^c$ can be made smaller than $\delta/2$ by choosing $m$ large.

Conditionally on $\mathcal{G}_n$, the difference in (121) is positive and it increases with $k$. It is therefore bounded by

$$S(j^n_{Th(n)^2}) - \sum_{j=0}^{Th(n)^2-1} s^n_j, \tag{122}$$

which is simply the time spent in $T^\varepsilon$ and $T_M$ during the first $j^n_{Th(n)^2}$ parts.

We first show that the time spent in $T^\varepsilon$ is small.

LEMMA 4.9. *For any $\delta > 0$ there exists $\varepsilon$ such that for a.e. $\boldsymbol{\tau}$ and $n$ large enough,*

$$\mathbb{P}\left[\left\{\sum_{i=0}^{j^n_{Th(n)^2}} e_i \tau_{Y(i)} \mathbb{1}\{Y(i) \in T^\varepsilon\} \geq 2^{n/\alpha}\delta\right\} \cap \mathcal{G}_n \Big| \boldsymbol{\tau}\right] \leq \delta. \tag{123}$$

PROOF. On $\mathcal{G}_n$ the first $Th(n)^2$ parts of the trajectory stays in $\mathbb{D}(n)$. The probability in (123) is thus bounded from above by

$$\mathbb{P}\left[\sum_{i=0}^{\zeta_n} e_i \tau_{Y(i)} \mathbb{1}\{Y(i) \in T^\varepsilon\} \geq 2^{n/\alpha}\delta \Big| \boldsymbol{\tau}\right], \tag{124}$$

where $\zeta_n$ is the exit time of $Y$ from $\mathbb{D}(n)$ [see (18)]. We show that there exists a constant $K_1$ independent of $\varepsilon$ such that for a.e. $\boldsymbol{\tau}$ and $n$ large enough

$$\mathbb{E}\left[\sum_{i=0}^{\zeta_n} e_i \tau_{Y(i)} \mathbb{1}\{Y(i) \in T^\varepsilon\} \Big| \boldsymbol{\tau}\right] \leq K_1 \varepsilon^{1-\alpha} 2^{n/\alpha}. \tag{125}$$

The claim of the lemma then follows by the Markov inequality.

To prove (125) we bound the expected time spent in traps with $\tau_x \leq 1$ first,

$$\mathbb{E}\left[\sum_{i=0}^{\zeta_n} e_i \tau_{Y(i)} \mathbb{1}\{\tau_{Y(i)} \leq 1\} \Big| \boldsymbol{\tau}\right] = \sum_{x \in \mathbb{D}} \tau_x G_{\mathbb{D}}(0,x) \mathbb{1}\{\tau_x \leq 1\}$$

$$\leq \sum_{x \in \mathbb{D}} G_{\mathbb{D}}(0,x) = \mathbb{E}(\zeta_n) = O(2^n) \ll 2^{n/\alpha}. \tag{126}$$

We divide the remaining part of $T^\varepsilon$ into disjoint sets $T^{\varepsilon 2^{-i+1}}_{\varepsilon 2^{-i}}$, where $i \in \{1,\ldots,i_{\max}\}$ and $i_{\max}$ is an integer satisfying

$$1/2 \leq 2^{-i_{\max}} \varepsilon g(n) < 1. \tag{127}$$

From condition (2) it can be showed easily that the probability that a fixed site $x$ is in $T^{\varepsilon 2^{-i+1}}_{\varepsilon 2^{-i}}$ is bounded by

$$p_{n,i} := \mathbb{P}[x \in T^{\varepsilon 2^{-i+1}}_{\varepsilon 2^{-i}}] \leq \mathbb{P}[\tau_x \geq 2^{-i}\varepsilon g(n)] \leq c\varepsilon^{-\alpha} g^{-\alpha} 2^{i\alpha}. \tag{128}$$

For any fixed $i \in \{1, \ldots, i_{\max}\}$ and $K'$ large we can write

$$
\mathbb{P}\left[\mathbb{E}\left[\sum_{j=0}^{\zeta_n - 1} e_j \tau_{Y(j)} \mathbb{1}\{Y(j) \in T_{\varepsilon 2^{-i}}^{\varepsilon 2^{-i+1}}\} \Big| \boldsymbol{\tau}\right] \geq K' \varepsilon^{1-\alpha} 2^{i(\alpha-1)} 2^{n/\alpha}\right]
$$

$$
(129) \qquad = \mathbb{P}\left[\sum_{x \in \mathbb{D}} G_{\mathbb{D}}(0, x) \tau_x \mathbb{1}\{x \in T_{\varepsilon 2^{-i}}^{\varepsilon 2^{-i+1}}\} \geq K' \varepsilon^{1-\alpha} 2^{i(\alpha-1)} 2^{n/\alpha}\right]
$$

$$
\leq \mathbb{P}\left[\sum_{x \in \mathbb{D}} G_{\mathbb{D}}(0, x) \mathbb{1}\{x \in T_{\varepsilon 2^{-i}}^{\varepsilon 2^{-i+1}}\} \geq K' \varepsilon^{-\alpha} 2^{i\alpha-1}\right].
$$

Using the Markov inequality (with $\lambda_n > 0$) this can be bounded by

$$
(130) \qquad \leq \exp(-\lambda_n K' \varepsilon^{-\alpha} 2^{i\alpha-1}) \prod_{x \in \mathbb{D}} [(1 - p_{n,i}) + p_{n,i} e^{\lambda_n G_{\mathbb{D}}(0,x)}]
$$

$$
\leq \exp(-\lambda_n K' \varepsilon^{-\alpha} 2^{i\alpha-1}) \prod_{x \in \mathbb{D}} [1 + c2^{i\alpha} g^{-\alpha} \varepsilon^{-\alpha} (e^{\lambda_n G_{\mathbb{D}}(0,x)} - 1)].
$$

Since $x \geq \log(1 + x)$, we have

$$
(131) \qquad \log \prod_{x \in \mathbb{D}} [1 + c2^{i\alpha} g^{-\alpha} \varepsilon^{-\alpha} (e^{\lambda_n G_{\mathbb{D}}(0,x)} - 1)]
$$

$$
\leq \sum_{x \in \mathbb{D}} c2^{i\alpha} g^{-\alpha} \varepsilon^{-\alpha} (e^{\lambda_n G_{\mathbb{D}}(0,x)} - 1).
$$

Let $\lambda_n = n/2G_{\mathbb{D}}(0,0)$. We divide the last sum into two parts. First, we sum over the sites that are close to the origin, $|x| \leq n^{2/(d-2)}$. Since $G_{\mathbb{D}}(0, x) \leq G_{\mathbb{D}}(0, 0)$, we have

$$
(132) \qquad \sum_{x \in D(n^{2/(d-2)})} c2^{i\alpha} g^{-\alpha} \varepsilon^{-\alpha} (e^{\lambda_n G_{\mathbb{D}}(0,x)} - 1)
$$

$$
\leq Cn^{2d/(d-2)} 2^{i\alpha-n} \varepsilon^{-\alpha} e^{\lambda_n G_{\mathbb{D}}(0,0)}
$$

$$
\leq Cn^{2d/(d-2)} 2^{i\alpha} 2^{-n} \varepsilon^{-\alpha} e^{n/2}.
$$

The last expression tends to 0 as $n \to \infty$.

By (146), $G_{\mathbb{D}}(0, x) \leq cn^{-2}$ for $x \in \mathbb{D}(n) \setminus D(n^{2/(d-2)})$. Therefore, the argument of the exponential in (131) is smaller than $c'n^{-1}$. Using the fact that $e^x - 1 \leq 2x$ for $x$ sufficiently close to 0 we get $e^{\lambda_n G_{\mathbb{D}}(0,x)} - 1 \leq cnG_{\mathbb{D}}(0,x)$ and thus

$$
(133) \qquad \sum_{x \in \mathbb{D} \setminus D(n^{2/(d-2)})} c2^{i\alpha} g^{-\alpha} \varepsilon^{-\alpha} (e^{\lambda_n G_{\mathbb{D}}(0,x)} - 1)
$$

$$
\leq \sum_{x \in \mathbb{D} \setminus D(n^{2/(d-2)})} Cn2^{i\alpha} \varepsilon^{-\alpha} g^{-\alpha} G_{\mathbb{D}}(0,x) \leq C2^{i\alpha} \varepsilon^{-\alpha} n.
$$

Here we used $\sum_{x \in \mathbb{D}} G_{\mathbb{D}}(0, x) = O(r(n)^2) = O(g^\alpha)$. From (132) and (133) it follows that the expression in (130) can be bounded from above by

$$(134) \qquad \exp(-K'cn\varepsilon^{-\alpha}2^{i\alpha})\exp(Cn\varepsilon^{-\alpha}2^{i\alpha}).$$

Therefore, it is possible to choose $K'$ large enough such that this bound decreases exponentially with $n$ for all $i \in \{0, \ldots, i_{\max}\}$.

Summing over all possible values of $i$ gives

$$
\begin{aligned}
(135) \qquad &\mathbb{P}\left[\bigcup_{i=0}^{i_{\max}}\left(\mathbb{E}\left[\sum_{i=0}^{\zeta_n-1} e_i \tau_{Y(i)} \mathbb{1}(Y(i) \in T_{\varepsilon2^{-i}}^{\varepsilon2^{-i+1}}) \Big| \boldsymbol{\tau}\right] \geq K'\varepsilon^{1-\alpha}2^{i(\alpha-1)}2^{n/\alpha}\right)\right] \\
&\leq ne^{-cn}.
\end{aligned}
$$

The Borel–Cantelli lemma then yields

$$(136) \qquad \mathbb{E}\left[\sum_{i=0}^{\zeta_n-1} e_i \tau_{Y(i)} \mathbb{1}(Y(i) \in T_{\varepsilon2^{-i}}^{\varepsilon2^{-i+1}}) \Big| \boldsymbol{\tau}\right] \leq K'\varepsilon^{1-\alpha}2^{i(\alpha-1)}2^{n/\alpha}$$

$\boldsymbol{\tau}$-a.s. for all $i$ and for $n$ large enough. Combining (126) and (136) we get easily (125). This completes the proof of Lemma 4.9. $\quad\square$

We show now that the set $T^M$ can be safely ignored.

Lemma 4.10. *For every $\delta$ there exist $m$ and $M$ such that for a.e. $\boldsymbol{\tau}$ and $n$ large enough*

$$(137) \qquad \mathbb{P}[\{Y \text{ hits } T_M \text{ before } j_{Th(n)^2}^n\} \cap \mathcal{G}_n | \boldsymbol{\tau}] \leq \delta.$$

Proof. As in the proof of Lemma 4.9 we can replace $j_{Th(n)^2}^n$ by $\zeta_n$. We use again the Borel–Cantelli lemma,

$$
\begin{aligned}
(138) \qquad &\mathbb{P}[\mathbb{P}[Y \text{ hits } T_M(n) \text{ before } \zeta_n | \boldsymbol{\tau}] \geq \delta] \\
&\leq e^{-\lambda_n \delta}\mathbb{E}[\exp\{\lambda_n \mathbb{P}[Y \text{ hits } T_M(n) \text{ before } \zeta_n | \boldsymbol{\tau}]\}].
\end{aligned}
$$

However,

$$
\begin{aligned}
(139) \qquad &\log \mathbb{E}[\exp\{\lambda_n \mathbb{P}[Y \text{ hits } T_M(n) \text{ before } \zeta_n | \boldsymbol{\tau}]\}] \\
&\leq \log \mathbb{E}\left[\exp\left\{\lambda_n \sum_{x \in \mathbb{D}} \mathbb{P}[Y \text{ hits } x \text{ before } \zeta_n]\mathbb{1}(x \in T_M)\right\}\right].
\end{aligned}
$$

Since $\mathbb{P}[x \in T_M] \leq cM^{-\alpha}g^{-\alpha}$, we get

$$
\begin{aligned}
(140) \qquad &\leq \sum_{x \in \mathbb{D}} \log\{1 + cM^{-\alpha}g^{-\alpha}(\exp\{\lambda_n \mathbb{P}[Y \text{ hits } x \text{ before } \zeta_n]\} - 1)\} \\
&\leq \sum_{x \in \mathbb{D}} cM^{-\alpha}\gamma^{-\alpha}\{\exp(\lambda_n \mathbb{P}[Y \text{ hits } x \text{ before } \zeta_n]) - 1\}.
\end{aligned}
$$

We choose $\lambda_n = n/2$ and divide the sum into two parts. For $|x| \leq n^{2/(d-2)}$ we use $\mathbb{P}[Y \text{ hits } x \text{ before } \zeta_n] \leq 1$. Hence,

$$
\begin{aligned}
(141) \quad & \sum_{x \in D(n^{2/(d-2)})} cM^{-\alpha}\gamma^{-\alpha}\{\exp(\lambda_n \mathbb{P}[Y \text{ hits } x \text{ before } \zeta_n]) - 1\} \\
& \leq cn^{2d/(d-2)}2^{-n}e^{n/2},
\end{aligned}
$$

which becomes negligible as $n \to \infty$.

By (148), for $|x| \geq n^{2/(d-2)}$ the argument of the exponential in (139) is smaller than $cn^{-1}$ and thus

$$
(142) \qquad \exp(\lambda_n \mathbb{P}[Y \text{ hits } x \text{ before } \zeta_n]) - 1 \leq cn|x|^{2-d}
$$

for some large $c$. We have thus

$$
\begin{aligned}
(143) \quad & \sum_{x \in \mathbb{D} \setminus D(n^{2/(d-2)})} cM^{-\alpha}g^{-\alpha}\{\exp(\lambda_n \mathbb{P}[Y \text{ hits } x \text{ before } \zeta_n]) - 1\} \\
& \leq cM^{-\alpha}g^{-\alpha}n \sum_{x \in \mathbb{D} \setminus D(n^{2/(d-2)})} |y|^{2-d} \leq cM^{-\alpha}n.
\end{aligned}
$$

Inserting (141) and (143) into (138) we get

$$
(144) \qquad \mathbb{P}[\mathbb{P}[Y \text{ hits } T_M(n) \text{ before } \zeta_n | \boldsymbol{\tau}] \geq \delta] \leq c\exp(-n\delta + c'M^{-\alpha}n).
$$

The proof is complete by taking $M$ large enough. $\quad\square$

## APPENDIX: PROPERTIES OF THE SIMPLE RANDOM WALK

We summarize here some useful facts about the Green's function and hitting probabilities of the simple random walk in the large ball $D(r) \subset \mathbb{Z}^d$, $d \geq 3$. The following lemma is taken from [11], Proposition 1.5.9.

LEMMA A.1.   *The Green's function $G_{D(r)}(\cdot, \cdot)$ of the simple random walk killed on exit from the ball $D(r)$ satisfies*

$$
(145) \qquad G_{D(r)}(0,0) = G_d(0) - O(r^{2-d})
$$

*and*

$$
(146) \qquad G_{D(r)}(0,x) = a_d(|x|^{2-d} - r^{2-d}) + O(|x|^{1-d}),
$$

*where $a_d = \frac{d}{2}\Gamma(\frac{d}{2} - 1)\pi^{-d/2}$.*

The hitting probabilities are controlled by the following lemma.

LEMMA A.2. *Let $p_r(0, x)$ denote the probability that the simple random walk started at 0 hits $x$ before exiting $D(r)$. The function $p_r(0, x)$ satisfies*

$$(147) \qquad p_r(0, x) \geq \frac{a_d}{G_d(0)}(|x|^{2-d} - r^{2-d}) + O(|x|^{1-d}),$$

$$(148) \qquad p_r(0, x) \leq a_d(|x|^{2-d} - r^{2-d}) + O(|x|^{1-d}).$$

*More precisely $p_n(0, x)$ can be bounded from above by*

$$(149) \quad p_r(0, x) \leq \frac{a_d}{G_d(0)}(|x|^{2-d} - r^{2-d} + O(|x|^{1-d}))(1 + O((n - |x|)^{2-d})).$$

PROOF. The first two claims follow from equation (146),

$$(150) \qquad G_{D(r)}(0, x) = p_r(0, x)G_{D(r)}(x, x),$$

and from $1 \leq G_{D(r)}(x, x) \leq G_d(0)$. The third fact is a consequence of (150) and

$$(151) \qquad \begin{aligned} G_{D(r)}(x, x)^{-1} &\leq G_{D_x(r-|x|)}(x, x)^{-1} = G_{D(r-|x|)}(0, 0)^{-1} \\ &= G_d(0)^{-1} + O((r - |x|)^{2-d}), \end{aligned}$$

which is a consequence of (145). $\quad\square$

**Acknowledgment.** We thank Takashi Kumagai for suggesting this problem of the scaling limit in $d \geq 2$.

COURANT INSTITUTE FOR MATHEMATICAL SCIENCES
NEW YORK UNIVERSITY
251 MERCER STREET
NEW YORK, NEW YORK 10012-1185
USA
E-MAIL: benarous@cims.nyu.edu

ÉCOLE POLYTECHNIQUE
FÉDÉRALE DE LAUSANNE
STATION 14 (CMOS)
1015 LAUSANNE
SWITZERLAND
E-MAIL: jiri.cerny@epfl.ch